\documentclass[10pt]{article}

\usepackage{amsfonts}
\usepackage{amsmath}
\usepackage{theorem}
\newtheorem{theorem}{Theorem}
\newtheorem{definition}[theorem]{Definition}
\usepackage[unicode,colorlinks,plainpages=false]{hyperref}
\newcommand{\openbox}{$\begin{array}{c}
\hspace*{-0.55em}\sqcap \hspace*{-0.60em}\\[-0.4em] \hline
\multicolumn{1}{c}{\hspace*{-0.60em}}\\[-0.8em]
\end{array}$}

\begin{document}

\centerline{\bf REFLEXIVE UNITARY SUBSEMIGROUPS OF}
\centerline{\bf LEFT SIMPLE SEMIGROUPS}

\bigskip

\bigskip

\centerline{A. Nagy}

\bigskip

\bigskip

\bigskip

\centerline{\bf Introduction}

\bigskip

Ideal series of semigroups \cite{1} play an important role in the examination of semigroups which have proper two-sided ideals. But the corresponding theorems cannot be used when left simple (or right simple or simple) semigroups are considered. So it is a natural idea that we try to use the group-theoretical methods (instead of the ring-theoretical ones) for the examination of these semigroups.

The purpose of this paper is to find a suitable type of subsemigroups of left simple semigroups which makes possible for us to generalize some notions (the notion of the normal series and the composition series of groups) and some results concerning the groups to the left simple semigroups. We note that the subsemigroups we are looking for are the reflexive unitary subsemigroups of left simple semigroups.

For notations and notions not defined here, we refer to \cite{1} and \cite {2}.

\bigskip

\bigskip

\centerline{\bf Reflexive unitary subsemigroups}

\bigskip

\bigskip

As it is known \cite{1}, a subset $H$ of a semigroup $S$ is said to be reflexive in $S$ if $ab\in H$ implies $ba\in H$ for all $a, b\in S$.

We say that a subset $U$ of a semigroup $S$ is left [right] unitary in $S$ (see \cite{1}) if $ab, a\in U$ implies $b\in U$ [$ab, b\in U$ implies $a\in U$] for all $a, b\in S$. A subset of $S$ which is both left and right unitary in $S$ is said to be unitary in $S$.

We note that if $A\subseteq B$ are subsemigroups of a semigroup $S$ and $B$ is unitary in $S$, then $A$ is unitary in $B$ if and only if $A$ is unitary in $S$. In this paper we shall use this fact without comment.

It can be easily verified that, in a group, the unitary subsemigroup are exactly the subgroups   , and the reflexive unitary subsemigroups are exactly the normal subgroups. So the notion of the reflexive unitary subsemigroup of a semigroup can be considered as a generalization of the notion of the normal  subgroup of a group.

\newpage

Reflexive unitary subsemigroups are important in the description of the group congruences of left simple semigroups. To show this importance, we need the following.

Let $S$ be a semigroup and $H$ a non-empty subset of $S$. Let
\[H... a=\{(s,t)\in S\times S:\quad sat\in H\}, \quad a\in S.\]
We can define a relation $P_H(S)$ on $S$ as follows:
\[P_H(S)=\{(a, b)\in S\times S:\quad H... a=H... b\}.\]
It can be easily verified that $P_H(S)$ is a congruence on $S$ such that $W^H=\{c\in S:\ H...c=\emptyset\}$ is a $P_H(S)$-class and an ideal of $S$. $P_H(S)$ is called the principal congruence on $S$ determined by $H$ (see \cite{1}).

If we consider and fix a semigroup $S$ then, for a subset $H$ of a subsemigroup $N$ of $S$, the principal congruence on $N$ determined by $H$ will be denoted by $P_H(N)$. we shall use $P_H$ instead of $P_H(S)$. For short, we shall use $N/P_H$ instead of $N/P_N(H)$. $P_H|N$ will denote the restriction of $P_H$ to $N$, that is, $P_H|N=P_H(S)\cap (N\times N)$.

We note that if $\alpha$ is a congruence on a semigroup $S$ then we shall not distinguish the congruence classes of $S$ modulo $\alpha$ from the elements of the factor semigroup $S/\alpha$.

Let $S$ denote a left simple semigroup. Then, by Theorem 10.34, Corollary 10.35 and Exercise for \S 10.2 of \cite{1}, if $H$ is a reflexive unitary subsemigroup of $S$ then the factor semigroup $S/P_H$ is a group with identity element $H$ and, conversely, if $P$ is a congruence on $S$ such that $S/P$ is a group with identity $H$ then $H$ is a reflexive unitary subsemigroup of $S$ and $P=P_H$.

The following theorem is important for the next.

\medskip

\begin{theorem}\label{th1} Every right unitary subsemigroup of a left simple semigroup is left simple.
\end{theorem}

{\it Proof.} Let $N$ be a right unitary subsemigroup of a left simple semigroup $S$. Consider two arbitrary elements $a$ an $b$ in $N$. Then there is an element $s$ in $S$ such that $sa=b$. As $N$ is right unitary in $S$, it follows that $s\in N$ which implies that $N$ is left simple. Thus the theorem is proved.\hfill\openbox

\medskip

The following theorem shows that, fixing a reflexive unitary subsemigroup $H$ of a left simple semigroup $S$, there is a one-to-one correspondence between the set of all unitary subsemigroups of $S$ containing $H$ and the set of all subgroups of $S/P_H$. In this connection the reflexive unitary subsemigroups correspond to the normal subgroups of $S/P_H$.

\medskip

\begin{theorem}\label{th2} Let $H\subseteq N$ be subsemigroups of a left simple semigroup $S$ such that $H$ is reflexive and unitary in $S$. Then $N$ is unitary in $S$ if and only if $N$ is saturated by $P_H$ and $N/P_H$ is a subgroup of $S/P_H$. In this case $P_H(N)=\break P_H|N$.
If $N$ is unitary in $S$, then it is also reflexive in $S$ if and only if $N/P_H$

\newpage
\noindent
is a normal subgroup of $S/P_H$. If this is the case, then $(S/P_H)/(N/P_H)$ is isomorphic with $S/P_N$.
\end{theorem}

{\it Proof}. Let $H$ be a reflexive unitary subsemigroup of a left simple semigroup $S$. Assume that $N$ is an unitary subsemigroup of $S$ with $H\subseteq N$. Let $a\in N$ be arbitrary. As $N$ is left simple by Theorem~\ref{th1}, there are elements $u, v\in N$ such that $uav\in H$. As $H$ is reflexive in $S$, $vua\in H$. Since $N$ is unitary in $S$, $vu\in N$. So, for every $b\notin N$, it follows that $vub\notin N$ from which we get $vub\notin H$. Thus $ubv\notin H$ and so $(a,b)\notin P_H$ which means that $N$ is saturated by $P_H$.
Since the inclusion $P_H|N\subseteq P_H(N)$ follows from the definition of $P_H|N$ and $P_H(N)$, it is sufficient to show that $P_H(N)\subseteq P_H|N$. To show this last inclusion, let $a$ and $b$ be arbitrary elements in $N$ with $(a,b)\in P_H(N)$. We prove $(a,b)\in P_H|N$. Assume, in an indirect way, that $(a,b)\notin P_H|N$. Then there are elements $x, y\in S$ such that either $xay\in H$ and $xby\notin H$ or $xay\notin H$ and $xby\in H$; we consider only the first case. Since $H$ is reflexive in $S$, it follows that $yxa\in H$ and $yxb\notin H$. As $a\in N$ and $N$ is unitary in $S$, we get $yx\in N$. Thus, for an arbitrary element $h$ in $H$, it follows that $yxah\in H$ and $yxbh\notin H$. Since $yx$ and $h$ are in $N$, we have
$(a,b)\notin P_H(N)$ which contradicts the assumption $(a,b)\in P_H(N)$. Consequently, $(a,b)\in P_H|N$, that is, $P_H(N)=P_H|N$.
It is evident that $N/P_H$ is a subsemigroup of $S/P_H$. Let $[b]\in N/P_H$ be arbitrary, where $[b]$ denotes the congruence class of $S$ modulo $P_H$ containing $b$. From $[b][b]^{-1}\subseteq H$ it follows that $[b]^{-1}\subseteq N/P_H$, because $N$ is unitary in $S$ and $[b]\subseteq N$. Thus $N/P_H$ is a subgroup of $S/P_H$.

Conversely, if a subsemigroup $N$ with $H\subseteq N$ is saturated by $P_H$ and $N/P_H$ is a subgroup of $S/P_H$, then $N$ is a unitary subsemigroup of $S$. Next, consider two arbitrary elements $a$ and $b$ in $S$, and let $\omega$ denote the canonical homomorphism of $S$ onto $S/P_H$. Then $ab\in N$ if and only if $a\omega b\omega \in N\omega$. So $N$ is reflexive in $S$ if and only if $N\omega$ is a normal subgroup of $S/P_H$.

To prove the isomorphism between $S/P_N$ and $(S/P_H)/(N/P_H)=G$, let $\beta$ denote the canonical homomorphism of $S/P_H$ onto $G$. Then
\[S'=\{ s\in S:\ s\omega \beta \ \hbox{is the identity of $G$}\}\]
equals $N$ and the kernel of $\omega \beta$ is $P_N$ (see Theorem 10.34 of \cite{1}). Thus $S/P_N$ is isomorphic with $G$. Thus the theorem is proved.\hfill\openbox

\medskip

\begin{theorem}\label{th3} If $H$ and $N$ are subsemigroups of an arbitrary semigroup $S$ such that $H$ is reflexive and unitary in $S$ and $H\cap N\neq \emptyset$, then $H\cap N$ is a reflexive unitary subsemigroup of $N$ and $\langle H, N\rangle /P_H$ is isomorphic with $N/P_{H\cap N}$.
\end{theorem}

{\it Proof.} Let $H$ an $N$ be subsemigroups of a semigroup $S$ such that $H$ is reflexive and unitary in $S$ and $H\cap N\neq \emptyset$. It can be easily proved that $H\cap N$ is a reflexive unitary subsemigroup of $N$. We may assume $S=\langle H, N\rangle $. If $[N]$ denotes the union of the $P_H$-classes $A$ of $S$ with $A\cap N\neq \emptyset$, then $S=[N]$. Let $(a, b)\in P_H$
for some $a, b\in N$. Then, for every $x, y\in N$, $xay\in H\cap N$ if and only if $xby\in H\cap N$, because $(a,b)\in P_H$, $a, b\in N$ and $N$ is a subsemigroup. So \break

\newpage
\noindent
$(a, b)\in P_{H\cap N}$, that is, $P_H|N\subseteq P_{H\cap N}(N)$. Next we show that $P_{H\cap N}(N)\subseteq P_H|N$. Let $a, b\in N$ with $(a, b)\in P_{H\cap H}(N)$. Assume, in an indirect way, that $(a, b)\notin P_H|N$. Then there are elements $x, y\in S$ such that, for example, $xay\in H$ and $xby\notin H$. Since $H$ is the identity element $S/P_H$ and $S=\langle H, N\rangle $, there are elements $u$ and $v$ in $N$ with $(x, u)\in P_H$ and $(y, v)\in P_H$. Using $(xay, uav)\in P_H$ and $(xby, ubv)\in P_H$, we have $uav\in H$ and $ubv\notin H$. Since $u,v\in N$, the last result contradicts $(a, b)\in P_{H\cap N}(N)$. Thus $(a,b)\in P_H|N$. Consequently $P_H|N=P_{H\cap N}(N)$. Nothing that, for any $x\in S$, there exists $u\in N$ such that $(x, u)\in P_H$, we can see that $\langle H, N\rangle /P_H$ is isomorphic with $N/P_{H\cap N}$.\hfill\openbox

\medskip

\begin{theorem}\label{th4} If $H$ and $N$ are unitary subsemigroups of a left simple semigroup $S$ such that $H$ is reflexive in $S$, then $H\cap N\neq \emptyset$, $\langle H, N\rangle$ is a unitary subsemigroup of $S$ and $\langle H, N \rangle=HN$. If $N$ is also reflexive in $S$, then $HN$ is reflexive in $S$. In this case $\langle H, N\rangle =HN=NH$.
\end{theorem}

{\it Proof.} Let $H$ and $N$ be unitary subsemigroups of a left simple semigroup $S$ such that $H$ is reflexive in $S$. Consider two elements $b_1$ and $b_2$ in $N$ with $(b_1, b_2)\in P_H$. By Theorem~\ref{th1}, $N$ is left simple. Thus there is an element $x$ in $N$ such that $xb_1=b_2$. So $(xb_1, ub_2)\in P_H$ for every $u\in H$. Since $P_H$ is a group congruence on $S$, we have $(x, u)\in P_H$ which means that $H\cap N\neq \emptyset$. Evidently $H\cap N$ is a unitary and reflexive subsemigroup of $N$. Let
\[K=\{ s\in S:\ (s, n)\in P_H\ \hbox{for some}\ n\in N\}.\]
Then $K$ is a subsemigroup of $S$. We show that $K$ is unitary in $S$. Let $a$, $b$ be arbitrary elements in $S$. Assume $a, ab\in K$. Then there is an element $a_1$ in $N$ such that $(a, a_1)\in P_H$. Since $H\cap N$ is a reflexive unitary subsemigroup of $N$ and $N$ is left simple, $P_{H\cap N}(N)$ is a group congruence on $N$. Thus there is an element $r$ in $N$ such that $ra_1\in H\cap N$. Then $ra\in H$ which implies $(rab, b)\in P_H$. Since $rab=r(ab)\in NK\subseteq K$, it follows that $b\in K$. Thus $K$ is left unitary in $S$. To show that $K$ is right unitary in $S$, assume $a, ba\in K$. Let $k$ denote the product $ba$. Then there are elements $a_2, k_1\in N$ such that $(a_2, a)\in P_H$ and $(k_1, k)\in P_H$. Since $N$ is left simple, there is an element $q$ in $N$ such that $qa_2=k_1$. Then $(ba, qa)\in P_H$. Since $P_H$ is a group congruence on $S$, it follows that $(b, q)\in P_H$, that is, $b\in K$. Thus $K$ is right unitary in $S$. So $K$ is unitary in $S$. The proof of the first assertion will be complete if we show that $HN=K$.
Since $H, N\subseteq K$, it follows that $HN\subseteq \langle H, N\rangle \subseteq K$. Let $k$ be an arbitrary element of $K$. Fix an element $v$ of $N$ such that $(v, k)\in P_H$. Since $S$ is left simple, there is an element $s$ in $S$ such that $sv=k$ which means that $s\in H$ and so $k\in HN$. Thus $K\subseteq HN$ and so $HN=\langle H, N\rangle =K$. Assume that $N$ is also reflexive in $S$. Then $H\cap N$ is a reflexive unitary subsemigroup of $S$ and both $H$ and $N$ are saturated by $P_{H\cap N}$. Denote $\omega$ the canonical homomorphism of $S$ onto $S/P_{H\cap N}$. Then $H\omega$ and $N\omega$ are normal subgroups of the group $S/P_{H\cap N}$ and so $N\omega H\omega=(NH)\omega$ is a normal subgroup of $S/P_{H\cap N}$.
Let $L=\{s\in S:$
\newpage
\noindent
$s\omega \in (NH)\omega\}$. Then, by Theorem~\ref{th2}, $L=NH=\langle H, N\rangle =HN$ is a reflexive unitary subsemigroup of $S$.\hfill\openbox

\medskip

\begin{theorem}\label{th5} (Zassenhaus lemma)  Let $S$ be a left simple semigroup and $A, B, M, N$ subsemigroups of $S$. Assume that $A$ and $B$ are unitary in $S$ and $A\cap B\neq \emptyset$. If $N$ and $M$ are reflexive unitary subsemigroups in $A$ and $B$, respectively, then $N(A\cap M)$ and
$M(B\cap N)$ are reflexive unitary subsemigroups in $N(A\cap B)$ and $M(A\cap B)$, respectively, such that $N(A\cap B)/P_{N(A\cap M)}$ is isomorphic with $M(A\cap B)/P_{M(B\cap N)}$.
\end{theorem}

{\it Proof.} First we prove that $A\cap M\neq \emptyset$. By the assumption for $A$ and $B$, $A\cap B$ is an unitary subsemigroup in $B$. Since $M$ is a reflexive unitary subsemigroup in $B$, we get $M\cap A=M\cap (A\cap B)\neq \emptyset$ (see also Theorem~\ref{th4}). Similarly, $B\cap N\neq \emptyset$. Using also Theorem~\ref{th3}, $A\cap M$ and $B\cap N$ are unitary subsemigroups of $A\cap B$. It can be easily verified that they are also reflexive in $A\cap B$. Then, by Theorem~\ref{th4}, $\langle A\cap M, B\cap N\rangle =(A\cap M)(B\cap N)=(B\cap N)(A\cap M)$ is a reflexive unitary subsemigroup in $A\cap B$. Let $C$, $D$ and $H$ denote the semigroups $A\cap B$, $(A\cap M)(B\cap N)$ and $C/P_D$. As $C$ is left simple by Theorem~\ref{th1}, $H$ is a group. We give a homomorphism $p$ of $NC$ onto $H$. Let $x$ be an arbitrary element of $NC$. Then there exist elements $a\in N$ and $c\in C$ such that $x=ac$. Let $xp$ be the congruence class of $C$ modulo $P_D(C)$ containing the element $c$. We prove that $p$ is well defined. Let $a_1\in N$ and $c_1\in C$ be arbitrary elements with $a_1c_1=x$. We prove that $(c, c_1)\in P_D(C)$. Since $H$ is a group whose identity element is $D$, there is an element $c'$ in $C$ such that $c_1c'\in D$. Thus $acc'=xc'=a_1c_1c'\in ND$. Since $D$ is unitary in $A\cap B$, it is unitary in $A$. By Theorem~\ref{th4}, $\langle H, D\rangle =ND$ is a unitary subsemigroup of $A$ and so of $S$. Since $acc'\in ND$ and $a\in \langle N, D\rangle =ND$ then it follows that $cc'\in ND$ and so $cc'\in ND\cap C$. As $ND\cap C=[((N\cap D)\cup (N\setminus D))D]\cap C=
[((N\cap D)D)\cap C]\cup [((N\setminus D)D)\cap C]= [((N\cap D)D)\cap C]\cup \emptyset =[((N\cap D)D)\cap C]=[(N\cap B)(N\cap B)(M\cap A)]\cap C=
[(N\cap B)(M\cap A)]\cap C=D\cap C=D$, we get $cc'\in D$. This result together with $c_1c'\in D$ imply that the $P_D(C)$-classes $[c]$ and $[c_1]$ of $C$ containing $c$ and $c_1$, respectively, the inverses of the $P_D(C)$-class $[c']$ of $C$ containing $c'$. As the factor semigroup $C/P_D$ is a group, we get $[c]=[c_1]$ which implies that $xp$ is well defined. Consequently, the mapping $p$ is well defined. To prove that $p$ is a homomorphism, let $x$ and $y$ be arbitrary element of $NC$. Then $xy\in NC$ because $NC=\langle N, C\rangle$ is a subsemigroup. Thus there are elements $a_1, a_2, a_3\in N\subseteq A$ and $c_1, c_2, c_3\in C\subseteq A$ such that $x=a_1c_1$, $y=a_2c_2$, $xy=a_3c_3$.
We show that $(c_1c_2, c_3)\in P_D(C)$. Since $a_1c_1a_2c_2=a_3c_3\in A$ then we have $(a_1c_1a_2c_2, a_3c_3)\in P_N(A)$. Since $N$ is the identity element of the group $A/P_N$ then it follows that $(c_1c_2, c_3)\in P_N(A)$. To prove $(c_1c_2, c_3)\in P_D(C)$, we have to that, for all
$u, v\in C=A\cap B$, both of the elements $uc_1c_2v$ and $uc_3v$ are either in $D$ or in $C\setminus D$. First we show that, for
\newpage
\noindent
every $t, s\in C$, the condition $(t, s)\in P_N(A)$ implies that both $t$ and $s$ are either in $D$ or in $C\setminus D$. Since
$N\subseteq \langle N  D\rangle =ND$ and $ND$ is a unitary subsemigroup of $A$ by Theorem 4, then we get that $ND$ is saturated by $P_N(A)$ by Theorem~\ref{th2}. Since $ND\cap C=D$ (see above) then $(t,s)\in P_N(A)$ implies either $t, s\in D$ or $t, s\notin D$ for every $t, s\in C$. Let $u$ and $v$ be arbitrary elements in $C$. As $u$ and $v$ re in $A$, taking into consideration that $(c_1c_2, c_3)\in P_N(A)$, we get $(uc_1c_2v, uc_3v)\in P_N(A)$. As $uc_1c_2v, uc_3v\in C$, both $uc_1c_2v$ and $uc_3v$ ar in either $D$ or in $C\setminus D$ (as above). Consequently $(c_1c_2, c_3)\in P_D(C)$ which means that $p$ is a homomorphism of $NC$ onto $H$. Recall that $D$ is the identity element of $H$. Let $Y=\{ y\in NC:\ yp=D\}$. We prove that $Y=N(A\cap M)$. Since $A\cap M\subseteq D$ (and so $Y\subseteq ND$), we have $N(A\cap M)\subseteq Y$. To prove $Y\subseteq N(A\cap M)$, let $y$ be an arbitrary element in $Y$. Then there are elements $a\in N$ and $d\in D$ with $y=ad$. Thus $y\in ND=N(B\cap N)(A\cap M)\subseteq N^2(A\cap M)=N(A\cap M)$, because $N$ is an unitary subsemigroup of $S$ and so (by Theorem~\ref{th1}) it is left simple, from which it follows that $N^2=N$. Thus $Y\subseteq N(A\cap M)$. Consequently $Y=N(A\cap M)$, that is, $N(A\cap M)$ is the identity element of the factor group $H\cong N(A\cap B)/ker_p$. Then $N(A\cap M)$ is a reflexive unitary subsemigroup of $N(A\cap B)$ and $N(A\cap B)/P_{N(A\cap M)}\cong H$. We can prove, in a similar way, that $M(B\cap N)$ is a reflexive unitary subsemigroup of $M(B\cap A)$ and $M(B\cap A)/P_{M(B\cap N)}\cong H$. By the transitivity of the isomorphism, the theorem is proved.\hfill\openbox

\bigskip

\bigskip

\centerline{\bf Normal series and composition series}

\bigskip

 \begin{definition}\label{df6} Let $S$ be a left simple semigroup. By a {\it normal series} of $S$ we mean a finite sequence
\begin{equation}\label{1}
S=S_0\supseteq S_1 \supseteq \dots \supseteq S_k
\end{equation}
of subsemigroups $S_i$ of $S$ with the property that $S_i$ is a reflexive unitary subsemigroup of $S_{i-1}$ for $i=1, \dots ,k$. The integer $k$ is called the length of the normal series (\ref{1}). The factor of (\ref{1}) are the groups $S_{i-1}/P_{S_i}$, $i=1, \dots , k$.

\noindent
We say that the normal series (\ref{1}) and a normal series
\begin{equation}\label{2}
S=H_0\supseteq H_1 \supseteq \dots \supseteq H_n
\end{equation}
are isomorphic if $n=k$ and there exists a permutation $i\mapsto i^*$ of the integers $i=1, \dots , k$ such that $S_{i-1}/P_{S_i}$ and $H_{i^*-1}/P_{H_{i^*}}$ are isomorphic for $i=1, \dots , k$.

\noindent
Moreover, we say that (\ref{2}) is a refinement of (\ref{1}) if $n\geq k$ and every $S_i$ equals some $H_j$, $i=1, \dots , k; j=1, \dots ,n$.
\end{definition}

\medskip

\begin{theorem}\label{th7} In a left simple semigroup $S$, every two normal series have refinements ending with the same subsemigroup of $S$.
\end{theorem}

\newpage

{\it Proof}. Let
\begin{equation}\label{3}
S=S_0\supseteq S_1 \supseteq \dots \supseteq S_k
\end{equation}
and
\begin{equation}\label{4}
S=H_0\supseteq H_1 \supseteq \dots \supseteq H_n
\end{equation}
be two normal series of a left simple semigroup $S$. Using Theorem~\ref{th4}, we get that $S_i\cap H_j\neq \emptyset$ for every $i=0, 1, \dots ,k$ and every $j=0, 1, \dots , n$. It can be easily verified that, for every $i=0, 1, \dots k$, $S_i\cap H_j$ is a reflexive unitary subsemigroup of $S_i\cap H_{j-1}$, $j=1, \dots ,n$. Similarly, for every $j=0, 1, \dots , n$, $S_i\cap H_j$ is a reflexive unitary subsemigroup of $S_{i-1}\cap H_j$, $i=1, \dots , k$. Thus
\begin{equation}\label{5}
S=S_0\supseteq S_1 \supseteq \dots \supseteq S_k \supseteq S_k\cap H_1 \supseteq S_k\cap H_2 \supseteq \dots \supseteq S_k\cap H_n
\end{equation}
and
\begin{equation}\label{6}
S=H_0\supseteq H_1 \supseteq \dots \supseteq H_n \supseteq H_n\cap S_1 \supseteq H_n\cap S_2 \supseteq \dots \supseteq H_n\cap S_k
\end{equation}
are two normal series of $S$ such that (\ref{5}) and (\ref{6}) are refinements of (\ref{3}) and (\ref{4}), respectively.\hfill\openbox

\medskip

\begin{theorem}\label{th8} (Schreier refinement theorem) In a left simple semigroup, every two normal series have isomorphic refinements.
\end{theorem}

{\it Proof}. Let
\begin{equation}\label{7}
S=S_0\supseteq S_1 \supseteq \dots \supseteq S_k
\end{equation}
and
\begin{equation}\label{8}
S=H_0\supseteq H_1 \supseteq \dots \supseteq H_n
\end{equation}
be two normal series of a left simple semigroup $S$. By Theorem~\ref{th7}, we may assume $S_k=H_n$. For every $i=1, \dots ,k$ and $j=0, 1, \dots , n$, let
\[S_{j,i}=S_i(S_{i-1}\cap H_j)=\langle S_i, (S_{i-1}\cap H_j)\rangle\]
and, for every $i=0, 1, \dots , k$ and $j=1, \dots , n$, let
\[H_{j,i}=H_j(H_{j-1}\cap S_i)=\langle H_j, (H_{j-1}\cap S_i)\rangle.\]
By Theorem~\ref{th4}, $S_{0,i}=S_i(S_{i-1}\cap H_0)=S_iS_{i-1}=\langle S_i, S_{i-1}\rangle =S_{i-1}$ and
$S_{n, i}=S_i(S_{i-1}\cap H_n)=S_i(S_{i-1}\cap S_k)=S_iS_k=\langle S_i, S_k\rangle =S_i$ for every $i=1, \dots , k$. Similarly, $H_{j,0}=H_{j-1}$ and $H_{j,k}=H_j$ for every $j=1, \dots , n$.
Theorem~\ref{th5} implies that, for every $i=1, \dots , k$ and $j=1, \dots , n$, the subsemigroups $S_{j,i}$ and $H_{j,i}$ are reflexive unitary subsemigroups of subsemigroups $S_{j-1,i}$ and $H_{j,i-1}$, respectively, such that $S_{j-1, i}/P_{S_{j,i}}$ is isomorphic with $H_{j, i-1}/P_{H_{ji}}$. Then
\begin{equation}\label{9}
S=S_{0,1}\supseteq S_{1,1}\supseteq \dots \supseteq S_{n,1}=S_1=S_{0,2}\supseteq S_{1,2}\supseteq \dots \supseteq S_{n,k}=S_k
\end{equation}
\newpage

\noindent
and
\begin{equation}\label{10}
S=H_{1,0}\supseteq H_{1,1}\supseteq \dots \supseteq H_{1,k}=H_1=H_{2,0}\supseteq H_{2,1}\supseteq \dots \supseteq H_{n,k}=H_n
\end{equation}
are isomorphic normal series of $S$ such that (\ref{9}) and (\ref{10}) are refinements of (\ref{7}) and (\ref{8}), respectively.\hfill\openbox

\medskip

\begin{definition}\label{df9} A normal series
$S=S_0\supseteq S_1 \supseteq \dots \supseteq S_k$ of a left simple semigroup $S$ will be called a composition series of $S$ if $S_{i-1}\neq S_i$ and $S_{i-1}$ has no reflexive unitary subsemigroup $T$ with $S_{i-1}\supset T\supset S_i$, $i=1, \dots , k+1$ (here $S_{k+1}$ denotes the empty set).
\end{definition}

\medskip

{\bf Remark}. It is clear that a normal series $S=S_0\supseteq S_1 \supseteq \dots \supseteq S_k$ of a left simple semigroup $S$ is a composition series
of $S$ if and only if $S_{i-1}\neq S_i$, every factor $S_{i-1}/P_{S_i}$ is a simple group ($i=1, \dots ,k$) and $S_k$ has no proper reflexive unitary subsemigroups.

\begin{theorem}\label{th10} (Jordan-H\"older theorem) If a left simple semigroup $S$ has a composition series, then every two composition series of $S$ are isomorphic with each other.
\end{theorem}

{\it Proof}. By Theorem~\ref{th8}, it is obvious.\hfill\openbox

\medskip

{\bf Example}. Let $(S_1;\circ )$, $(S_2; *)$ be left simple semigroups such that $S_1\cap S_2=\emptyset$ and there is an isomorphism $\alpha$ of the semigroup $S_1$ onto $S_2$. On the set $F=S_1\cup S_2$, define an operation as follows: for every $e, f\in F$, let
\[ef=\begin{cases} e\circ f, & \text{if $e, f\in S_1$}\\
e*f\alpha & \text{if $e\in S_2, f\in S_1$}\\
e\alpha *f & \text{if $e\in S_1, f\in S_2$}\\
a\alpha ^{-1}\circ f\alpha ^{-1} & \text{if $e, f\in S_2$.}
\end{cases}\]
It can be shown that $F$ is a left simple semigroup such that $S_1$ is a reflexive unitary subsemigroup of $F$.

\medskip

{\bf Remarks}. If
\begin{equation}\label{11}
S=S_0\supseteq S_1 \supseteq \dots \supseteq S_k
\end{equation}
is a normal series of a left simple semigroup $S$, then every $S_i$, $i=0, 1, \dots , k$ is a unitary subsemigroup of $S$ and so they are left simple semigroups (see Theorem~\ref{th1}). If (\ref{11}) is a composition series, then $S_k$ has no proper reflexive unitary subsemigroups. So it is a natural problem to describe the structure of left simple semigroups which have no proper reflexive unitary subsemigroups.

By Theorem 1.27 of \cite{1}, a semigroup is left simple and contains an idempotent element if and only if it is a direct product of a left zero semigroup and a group. So we can prove that a left simple semigroup with an idempotent element has no proper reflexive unitary subsemigroups if and only if it is a left zero semigroup.

\medskip

{\bf A problem}. Find all left simple semigroups without idempotent elements which have no proper reflexive unitary subsemigroups.

\bigskip

\noindent
Attila Nagy

\noindent
Department of Algebra

\noindent
Mathematical Institute

\noindent
Budapest University of Technology and Economics

\noindent
e-mail: nagyat@math.bme.hu

\end{document}